\documentclass{article}

\usepackage{amssymb, amsmath}
\usepackage{amsfonts, latexsym}
\usepackage{a4wide} 
\usepackage{scalerel}

\newtheorem{theorem}{Theorem}[section] 
\newtheorem{lemma}[theorem]{Lemma}     
\newtheorem{corollary}[theorem]{Corollary}
\newtheorem{prop}[theorem]{Proposition}

\newtheorem{definition}{Definition}

\newtheorem{example}{Example}

\newcommand{\stern}[1]{{#1}^*}

\newcommand\reallywidehat[1]{%
\begin{array}{c}
\stretchto{
  \scaleto{
    \scalerel*[\widthof{#1}]{\bigwedge}
    {\rule[-\textheight/2]{1ex}{\textheight}} 
  }{1.25\textheight} 
}{0.5ex}\\           
#1\\                 
\rule{-1ex}{0ex}
 \end{array}
 }

\newcommand{\G}{\mathcal{G}}
\newcommand{\F}{\mathcal{F}}
\newcommand{\M}{\mathcal{M}}
\newcommand{\Hh}{\mathcal{H}}
\newcommand{\comp}{\mbox{\^{o}}}

\title{Group construction in non-trivial geometric $C$-minimal structures.}
\author{F. Delon and F. Maalouf}
\date{}
\begin{document}
\maketitle

\begin{abstract}We show that an infinite group is definable in any non trivial geometric $C$-minimal structure which is definably maximal and does not have any definable bijection between a bounded interval and an unbounded one in its canonical tree. No kind of linearity is assumed.
\end{abstract}
\section{Introduction}

In the spirit of the construction of a field from its projective plane, Boris Zilber proposed an ambitious program: in model theory, the notion of algebraic closure in suitable first order structures gives rise to combinatorial geometries. These geometries can be very similar to projective spaces. Zilber conjectured that a strongly minimal structure interprets an infinite group, or even an infinite field, as soon as it fulfills some conditions, conditions that are clearly necessary \cite{Z}. This conjecture turned out to be false in general. However, together with Ehud Hrushovski they were able to establish  that the conjecture holds for what they called ``Zariski structures'', first order structures with a topology which mimics the Zariski topology \cite{HZ}. Ya'acov Peterzil and Sergei Starchenko proved a variant of the conjecture for the class of o-minimal structures \cite{PS}. O-minimal structures are linearly ordered structures, thus endowed with the topology defined by the ordering, and they present 
strong analogies with strongly minimal structures.
\par\ \par \noindent
It is then natural to ask the question for $C$-minimal structures. The $C$-minimality condition is  an equivalent of strong minimality in the setting of ultrametric structures (or more generally $C$-structures) just as o-minimality is an equivalent of strong minimality in the setting of ordered structures. However, the Steinitz exchange property which is a consequence of strong minimality and o-minimality, does not hold for all $C$-minimal structures. 
If we assume it, we are in the setting of \emph{geometric structures} as defined by Ehud Hrushovski and Anand Pillay \cite{HP}, since $C$-minimal structures do eliminate the quantifier $\exists^\infty$, the other required property. 
Geometric structures offer a common framework for strongly minimal or o-minimal structures as well as for many classical mathematical structures and provide tuned tools and techniques. 
\par\ \par\noindent
The second author had constructed an infinite definable group in any geometric $C$-minimal structure, which is non-trivial and locally modular (\cite{maalouf1} and \cite{maalouf2}). 
In this paper, we remove the assumption of local modularity. New arguments have to be blown into the proof and we follow the spirit of \cite{PS}. 
Some extra conditions are assumed, that appear in \cite{annalestoulouse} and in the context of fields in \cite{HK}. They in particular guarantee the existence of limits of unary functions in the neighborhood of a point.  This allows us to copy an essential element of \cite{PS}: the notion of ``tangent''. For a definable curve $X$ and a definable family $\cal F$ of curves, where $X$ and all the curves of $\cal F$ pass through some fixed point $P$, the idea is to determine a curve in $\cal F$ which, on a neighborhood of $P$, is closer to $X$ than any other element of $\cal F$.
\par\ \par \noindent
More precisely we show the following.\par\ \par\noindent
{\bf Theorem:}\emph{
   Let $\M$ be a $C$-minimal structure which is definably maximal, geometric and non-trivial. Suppose moreover that in the underlying tree of $\M$ there is no definable bijection between a bounded interval and an unbounded one. Then there is an infinite group definable in $\M$.}\par\ 
%
%
%
%
%
\par\ \noindent
The paper is organized as follows. 
Section 2 is devoted to preliminaries on $C$-minimal structures. 
In Section \ref{snicefamilies}, we show that either an infinite group or a family $\F$ of functions with \textit{nice properties} is definable in $\M$. The idea is then to get a group law by composing the elements of this family of functions. As the composition is in general not in $\F$, we have to \emph{approximate} it with a function from $\F$. For this end, we introduce in Sections \ref{stangents} and \ref{sderivatives} the notions of \emph{tangents} and \emph{derivatives}, and study their properties. We use them in Section \ref{sgrouplaw} to construct an infinite group in $\M$.

\section{Preliminaries}
We start with a few definitions and preliminary results. $C$-structures and $C$-minimal structures have been introduced and studied in  \cite{macphersonhaskell} and  \cite{macphersonsteinhorn}. We remind in what follows their definition and  principal properties.\par\ \par\noindent
{\bf\emph{Notations:}} We use $\mathcal{M}, \mathcal{N},...$ to denote structures and  $M,N,...$ for their underlying sets. \par\ \par\noindent
In this paper, a $C$-structure is a  structure $\mathcal{M}=(M,C,...)$, where $C$ is a ternary predicate satisfying the following axioms:
\begin{itemize}
\item $\forall x,y,z, C(x,y,z)\longrightarrow  C(x,z,y)$ 
\item $\forall x,y,z, C(x,y,z)\longrightarrow \neg C(y,x,z)$
\item $\forall x,y,z,w, C(x,y,z)\longrightarrow  [ C(x,w,z)\vee C(w,y,z) ] $
\item $\forall x,y, x\neq y,\, \exists z\neq y, C(x,y,z)$.
\end{itemize}
Note that these C-structures are sometimes called ``dense'', see  \cite{delonlms}. 
Let $\mathcal{M}$ be a $C$-structure. We call \emph{cone}\footnote{Note: $M$ is not a cone.} or \emph{open ball} any subset of $M$ of the form $\{x; \mathcal{M}\models C(a,x,b)\}$, where $a$ and $b$ are two distinct elements of  $\mathcal{M}$. We call $0$-\emph{level set} or \emph{closed ball} any subset of $M$ of the form $\{x; \mathcal{M}\models \neg C(x,a,b)\}$, where $a$ and $b$ are two elements of  $\mathcal{M}$. A set is said to be a \emph{ball} if it is an open ball or a closed ball. It follows from the first three axioms of $C$-relations that the cones of $M$ form a basis of a completely disconnected topology on $M$. The last axiom guarantees that all cones are infinite. \par\ \par\noindent
Let $(T,\leq)$ be a partially ordered set. We say that $(T,\leq)$ is a \emph{tree} if the set of elements of $T$ less than any fixed element is totally ordered by $\leq$, and if any two elements of $T$ have a greatest lower bound. A \emph{leaf} is a maximal element of $T$. A \emph{branch} is a maximal totally ordered subset of $T$. It is easy to check that if $a$ and $b$ are two branches of $T$, then  $sup(a\bigcap b)$ exists. On the set of branches of $T$, we define a ternary relation $C$ in the following way: we say that $C(a,b,c)$ is true if and only if  $\mathrm{sup}(a\bigcap b)<\mathrm{sup}(b\bigcap c)$. It is easy to check that this relation on the set of branches satisfies the first three axioms of a $C$-relation.\par\ \par\noindent
A theorem from \cite{adeleke} says that $C$-structures can be looked at as sets of branches of a tree, equipped with the $C$-relation as defined above. The construction of the \emph{underlying tree} $T(\mathcal{M})$ of a $C$-structure $\M$ has been slightly modified in \cite{delonlms}. In this new construction, $M$ can be identified with the set of leaves of $T(\mathcal{M})$. The tree $T(\M)$ appears as the quotient of $M^2$ by an adequate equivalence relation which is definable in $(M,C)$. To an element $x \in M$, we associate the branch 
$br_x := \{ \nu \in T(\mathcal{M}) : \nu \leq x \}$. 
To elements $x,y\in M$, we associate the node $t:=sup (x\cap y)$, where $x$ and $y$ are seen as branches of $T(\mathcal{M})$. This operation is well defined, and we say then that $x$ and $y$ \emph{branch at} $t$. If $a$ and $b$ are two distinct elements of $M$ branching at a node $\nu$, we denote  $\Lambda_\nu(b):=\{x\in M; C(a,x,b)\}$, we call it \emph{the cone of $b$ at $\nu$}, and we say that $\nu$ is its \emph{basis}. For a node or a leaf $\nu$ of $T(\mathcal{M})$, 
we denote by $\Lambda_\nu$ the 
closed ball of $M$ defined by $\nu$. This corresponds to the set of all elements of $M$ which contain $\nu$ when 
considered 
as branches of $T(\mathcal{M})$. We call $\nu$ \emph{the basis of} $\Lambda_\nu$.\par\ \par\noindent
For $a,b\in T(\M) \cup \{-\infty\}$ with $a<b$, the interval $(a,b)$ is said to be \emph{bounded from above} (respectively, \emph{bounded from below}) if $b$ is not a leaf (respectively, if $a\neq-\infty$).

\begin{definition}Let  $\M=(M,C,...)$ be a $C$-structure.
\begin{enumerate}
\item  The structure $\M$ is \emph{geometric} if for any structure $\mathcal{N}$ elementarily equivalent to $\M$, the algebraic closure in $\mathcal{N}$ has the exchange property. In this case, the algebraic closure is the closure operator of a pregeometry on $\mathcal{N}$. If this pregeometry is trivial on any $\mathcal{N}$ (i.e. $\mathrm{acl}\: A = \bigcup \{\mathrm{ acl}\:\! a ; a \mbox{ a singleton in } A\}$ for any $A \subset N$), then $\M$ is said to be \emph{trivial}. See \cite{HP} and \cite{Marker}. 
\item The structure $\M$ is \emph{definably maximal} if any definable family of cones which is linearly ordered by the inclusion has a non-empty intersection.
\item   The structure  $\M$ is $C$-\emph{minimal}  if and only if for any structure $\mathcal{N}=(N,C,...)$ elementarily equivalent to $\M$, any definable subset of  $N$ can be defined without quantifiers using  only the relations $C$ and $=$.
\end{enumerate}
\end{definition}

Two comments on these definitions.\\
1. Geometric structures are provided with notions of independence, dimension, and generic points. 
In other respects $C$-minimal structures admit a cellular decompositions (see \cite{macphersonhaskell} and its complement \cite{Pablo}), which gives rise to a topological dimension. In geometric $C$-minimal structures, these two dimensions coincide. This means that 
a point is generic exactly when any definable set containing it contains a box of the ambient space. \\
2. Let us define a linearly ordered structure to be definably maximal if any definable decreasing family of bounded closed intervals has a non-empty intersection. Then any o-minimal structure is definably maximal. But not any $C$-minimal structure is definably maximal \cite{annalestoulouse}, nor geometric \cite{macphersonhaskell}. 

\begin{prop}\label{nobadfunc}
Let $\M$ be a geometric  $C$-minimal structure and $T$ its underlying tree. Let $f: M\longrightarrow T\setminus M$ be a definable partial function. Then $\mathrm{dom}(f)$ can be written as a definable union $F\cup K$ such that $F$ is finite, and $f$ is locally constant on $K$.
\end{prop}
\emph{Proof:}
This is a direct consequence of Propositions 3.9 and 6.1 of \cite{macphersonhaskell}.
\hfill$\square$\par\ \par\noindent

\begin{lemma}\label{formulevoisinage}
Let $\M$ be a geometric $C$-minimal structure, $\varphi$ a formula in two variables over $\M$ and $b\in M$ generic over the parameters defining $\varphi$. Let $D$ be a cone containing $b$ and suppose that for all $u\in D$ there is a subcone $V_u$ of $D$ containing $u$ such that for all $v\in V_u$, $\M\models \varphi(u,v)$. Then there exists a neighborhood $D'$ of $b$ such that   for all $u,v\in D'$, $\M\models \varphi(u,v)$.
\end{lemma}
\emph{Proof:}
For $u\in D$, let $f(u)$ be the node on the branch $br_u$ of $u$ such that 
\[f(u)=\min \{ \mathrm{inf}\{\nu\in br_u; \forall v\in\Lambda_{\nu}:\M\models\varphi(u,v)\},\mbox{basis of }D\}.\] 
 The function $f: M \rightarrow T(\cal M)$ is definable and  
 the hypotheses on $\varphi$ imply $f(D) \subset T({\cal M}) \setminus M$, so by Proposition \ref{nobadfunc}, $\mathrm{dom}(f)$ can be written as a definable union $F\cup K$ such that $F$ is finite, and $f$ is locally constant on $K$. By genericity, $b\notin F$. Thus there is a cone $D'\subset D$ containing $b$ on which $f$ is constant. Therefore, for any $u,v\in D'$, the formula $\varphi(u,v)$ is satisfied in $\M$. 
 \hfill$\square$\par\ \par\noindent

\begin{lemma}\label{fonctionbornee}
Let $\M$ be a $C$-minimal structure, $c\in M$ and $\varphi(i,x)$ a formula in two variables on $\M$ where $i$ ranges in the sort $T(\M)$ in an interval $I = ]\rho,c[$ of the branch of $c$, and $x$ ranges in a cone $U\subset M$. We suppose that, for any $x,y\in M$, there is no definable bijection between a bounded interval of the branch of $x$ and an unbounded interval of the branch of $y$. Suppose that 
\[\forall x \in U, \exists i \in I, \forall j \in I \ [j\geq i \rightarrow \varphi(j,x)].\] Then 
\[\exists i \in I, \forall x \in U, \forall j \in I \ [j\geq i \rightarrow \varphi(j,x)].\]   
\end{lemma}
\emph{Proof:}
For $i\in I$, let $U_i$ be the subset of elements $x\in U$ such that $\forall j \in I \ [j \geq i \rightarrow  \varphi(j,x)]$. The map $i\mapsto U_i$ is increasing, and the union of all the $U_i$ contains $U$. Then for every ball $\Lambda_{\nu}\subsetneq U$, $\Lambda_{\nu}$ is strictly contained in the union of the $U_i$, which is an increasing and definable family of definable subsets of $M$. By $C$-minimality the $U_i$ are (uniformly definable) Swiss Cheese and $\Lambda_{\nu}$ must be contained in one of them. For every $\nu$ let $i_{\nu}$ be the infimum of the $i \in br_c$ with this property. Fix a node $\nu_0$ such that $\Lambda_{\nu_0}\subset U$ and consider the function $f$ which to a node between the basis of $U$ and $\nu_0$ associates the element $i_{\nu}$. 
The branches equipped with the structure induced by $\M$ are o-minimal, 
thus there is a finite definable partition of the domain of $f$ such that on each piece $f$ is either constant or bijective monotonous. Consequently, 
since the domain of $f$ is bounded, its image must be bounded from above by some $i_0\in br_c$. By the choice of $i_0$ we have $\forall x \in U, \forall j \in I \ [j\geq i_0 \rightarrow \varphi(j,x)]$.      
\hfill$\square$\par\ \par\noindent

\begin{corollary}\label{fonctionbornee2}
Let $\M$ be a $C$-minimal structure with no definable bijection between a bounded interval and an unbounded interval of $T(\M)$. 
Let $\{ D(\mu) : \mu \in I \}$ be a family of uniformly definable subsets of $M$, indexed by an interval $I$ of $T(\cal M)$, $I=]\rho,c[$ for some $c \in M, \rho \in br_c$, and such that\\
- $\nu < \mu$ implies $D(\nu) \subset D(\mu)$\\
- $\bigcup_{\nu \in I} D(\nu)$ is a cone $\Gamma$ in $M$. \\
Then $D(\nu_0)=\Gamma$ for some $\nu_0 \in I$. 
\end{corollary}


\section{Nice families of functions.}\label{snicefamilies}
Let $\M$ be a $C$-minimal structure. 
\begin{lemma}\label{CCC}
Let $V$ be a cone in $M$, $e$ an element of $V$, and $f,g,h:V \rightarrow M$ definable functions.   Then there is a neighborhood $W$ of $e$ such that either any element $x\in W\setminus\{e\}$ satisfies $C(f(x),g(x),h(x))$, or any element $x \in W\setminus\{e\}$ satisfies $\neg C(f(x),g(x),h(x))$. 
\end{lemma}
\emph{Proof:} 
Let $D\subset M$ be the set of elements $x$ of $V$ such that $C(f(x),g(x),h(x))$ holds. Either  $e$ is an accumulation point of  $D$, in which case by $C$-minimality there is a neighborhood  $W$ of $e$ such that $W\setminus\{e\}\subset  D$, or $e$ is an accumulation point of the complement $\neg D$ of $D$, in which case there is a neighborhood $W$ of $e$ such that $W\setminus\{e\}\subset \neg D$.\hfill$\square$\par\ \par\noindent
 {\bf\emph{Notation:}}\begin{enumerate}\item Let $x,y,z\in M$. We denote by $\Delta(x,y,z)$ the property \[\neg C(x,y,z)\wedge \neg C(y,x,z).\]
\item Let $V\subset M$ be a cone, $e$ an element of $V$, and $f, g, h$ functions from $V$ to $M$. We denote by $C_e(f,g,h)$ (or $C(f,g,h)$ if there is no confusion on $e$) the following property: there exists a neighborhood $W$ of $e$ such that \[\forall x\in W\setminus\{e\} : C(f(x),g(x),h(x)).\] We define $(\neg C)_e(f,g,h)$ and $\Delta_e(f,g,h)$ in the same way, and we denote them by $(\neg C)(f,g,h)$ and $\Delta(f,g,h)$ respectively, if there is no confusion on $e$.
                     \end{enumerate} 
\begin{definition} Let $V$ be a cone in $\M$, $e\in V$, and $f:V \rightarrow M$ a definable function.
\begin{enumerate}
 \item The function $f$ is said to be \emph{dilating on a neighborhood of $e$}, or just \emph{dilating} if there is no confusion, if it satisfies $ C_e(f\circ f,f,id_V)$.
\item The function $f$ is said to be \emph{non-dilating on a neighborhood of $e$}, or just \emph{non-dilating} if there is no confusion, if it satisfies $( \neg C)_e(f\circ f,f,id_V)$.
\end{enumerate}
\end{definition}
\begin{definition}\label{nice} Let $V \subset M$ be a cone and $\F=\{f_u:u\in U\}$ a definable family of definable functions from $V$ to $M$,  indexed by a cone $U\subset M$. The family $\F$ is said to be \emph{a nice family of functions} if there is an element $e\in V$  with the following properties: \begin{enumerate}
\item All the $f_u$ are $C$-automorphisms of the cone $V$.
  \item For every $u\in U$, we have $f_u(e)=e$.
 \item For any fixed $x\in V\setminus\{e\}$, the application $U\longrightarrow V$ which to $u$ associates $f_u(x)$, is a $C$-isomorphism from $U$ onto some subcone of $V$.      \hfill$(*)$                                                                                                                                                                                                                            \end{enumerate}
A nice family $\F$ of functions is said to have an \emph{identity element} if for some $u_0\in U$, $f_{u_0}=id_V$. 
\end{definition}
{\bf\emph{Notation:}} If $\F$ is a nice family of functions, possibly with identity, then $U,V, e$ and $u_0$ will be as in Definition \ref{nice} if there is no other precision.  For a set $W$ containing $e$, the set $W\setminus\{e\}$ will be denoted by $\stern{W}$.\par\ \par\noindent
{\bf\emph{Terminology:}} For a nice family of functions $\F$, the sets $U$ and $V$ are called respectively  the \emph{index set} and the \emph{domain}, and $e$ is called the \emph{ absorbing element}.

\begin{prop}\label{existencefamille}
Let $\M$ be a non trivial geometric $C$-minimal structure. Suppose furthermore that  in the underlying tree of $\M$ there is no definable bijection between a bounded interval and an unbounded one.
Then either an infinite group or a nice family of non-dilating functions with identity is definable in $\M$.

\end{prop}
\emph{Proof:}
Without loss of generality, we can suppose that $\M$ is $\aleph_1$-saturated. By Proposition 15 and Lemma 21 of \cite{maalouf1}, we can find cones $U_1,V_1$ in $M$ and a definable family of functions $\Hh=\{h_u:u\in U_1\}$ of (continuous) $C$-automorphisms of $V_1$ satisfying $(**)$: for every $x\in V_1$ the map  $v\mapsto h_v(x)$ is a continuous $C$-isomorphism from $U_1$ onto some cone. Fix a generic triplet $(a,b,e)\in U_1\times U_1\times V_1$, and let $e':=h_a\circ h_b(e)$. By $(**)$, for any $v \in U_1$ there is at most one $z \in U_1$ such that $h_z \circ h_v(e)=e'$. When such a $z$ exists, define $\theta(v) := z$. 
By $(**)$ again and the genericity of $b$, the function $\theta$ is well defined on some cone $U_2$ containing $b$. 
 \par\ \noindent \par

 Assume first that for all neighborhoods $U_a,U_b$ of $a,b$ there are $u \in U_a, u \not= a$  and $v \in U_b, v \not= b$ such that $h_u \circ h_v$ and $h_a \circ h_b$ coincide on  some neighborhood of $e$ as soon as they agree on $e$. \\[5 mm]
 Claim. Under this assumption there are cones $U_a$, $U_b$ and $V$ containing  $a,b$ and $e$ respectively, such that for all $(u,v),(u',v') \in U_a \times U_b$, 
 $h_u\circ h_v$ and $h_{u'} \circ h_{v'}$ coincide on $V$ as soon they agree on $e$. \\[2 mm]
Proof. Take a cone $U_3\subset U_1$ containing $b$ such that all  $v \in U_3\setminus \{b\}$ have the same type on $(a,b,e)$.  Necessarily $U_3 \subset U_2$. This implies that, for any $v\in U_3$ there is $u \in U_1$ such that  
$h_u \circ h_v (e) = e'$, but then $u=\theta(v)$.  
So it follows from our assumption that for all $v\in U_3$, $h_{\theta(v)} \circ h_v$ and $h_a \circ h_b$ have the same germ on $e$.
Define $U_b' = U_3$ and $U_a' = \theta(U_b')$. For any $(u,v) \in U_a' \times U_b'$ there is some neighborhood $W$ of $e$ such that  $h_u \circ h_v$ and $h_a \circ h_b$ coincide on $W$ as soon as they agree on $e$. 
\\ 
Let $\varphi$ be the formula defined as follows: $\varphi (u,v,\nu) :\longleftrightarrow$ ``$h_u \circ h_v$ and $h_a \circ h_b$ coincide on $\Lambda_\nu(e)$ as soon as they agree on $e$''. We first fix $u\in U'_a$ and apply Lemma \ref{fonctionbornee} to get $\nu_u$ such that for all $v\in U'_b$, $h_u \circ h_v$ and $h_a \circ h_b$ coincide on $\Lambda_{\nu_u}(e)$ as soon as they agree on $e$. We apply Lemma \ref{fonctionbornee} again to get a cone  $V$ containing $e$ which satisfies the following:  for all $(u,v)\in U_a'\times U_b'$, $h_u \circ h_v$ and $h_a \circ h_b$ coincide on $V$ as soon as they agree on $e$. 
We can define $V$ with parameters, say $c$, independent over $(a,b,e)$. 
Thus there are cones $U_a$ and $U_b$  containing respectively $a$ and $b$ such that any point in $U_a \times U_b$ has same type as $(a,b)$ over $(c,e)$. Therefore, for all $(u,v),(u',v') \in U_a \times U_b$, 
 $h_u\circ h_v$ and $h_{u'} \circ h_{v'}$ coincide on $V$ as soon they agree on $e$.
\hfill $\dashv$
\par\ \par\noindent
In this first case, we construct first an infinite $C$-group $G$ type-definable in $\M$. For this end, we proceed exactly as in the proof of Theorem 19 of \cite{maalouf1}, applying the above claim instead of Lemma 22 of \cite{maalouf1}. Indeed the local modularity was only used for proving Lemma 22. Then Theorem 1 of \cite{maalouf2} gives us an infinite subgroup of $G$ definable in $\M$.

 \par\ \par 
  Assume now that we are not in the above case. By $C$-minimality, there are cones  $U_a,U_b$  containing respectively $a,b$ such that, for all $u \in U_a, u \not= a$ and $v \in U_b, v \not= b$ with $h_u \circ h_v(e)=e'$, the graphs of $h_u \circ h_v$ and $h_a \circ h_b$ do not intersect on $\stern{W}$ for some cone $W \ni e$. We can suppose without loss of generality that $U_b\subset U_2$. For $u\in U_2$, define $g_u:=h_{\theta(u)}\circ h_u$.
 \par\
 The element $b$ has the following property: there is a cone $U \ni b$ (for example $U=U_b$) such that for all $u\in U\setminus\{b\}$, the graphs of $g_u$ and $g_b$ do not intersect on $\stern{W}$ for some cone $W \ni e$. By genericity, this property holds on some cone $U_3\subset U_b$ containing $b$. Hence we have the following: for all $v\in U_3$, there is a cone $U_v \ni v$ such that for all $u\in U_v\setminus\{v\}$, the graphs of $g_u$ and $g_v$ do not intersect on $\stern{W}$ for some cone $W \ni e$. By Lemma \ref{formulevoisinage}, there is some cone $U_4\subset U_3$ containing $b$ such that 
 for any $u, v\in U_4$, $u\neq v$, there is some $W \ni e$ such that the functions $g_{u}$ and $g_{v}$ agree on exactly one point of $W$, namely $e$. \par\
For $\nu \in br_e$ and $v\in U_4$, we define \[D(\nu,v) := \{v\}\cup\{ u \in U_4 : \mathrm{Gr}(g_v) \cap \mathrm{Gr}(g_u) \cap [\Lambda_\nu(e)^*\times M] = \emptyset \},\] 
where $\mathrm{Gr}(g_u)$ and $\mathrm{Gr}(g_v)$ are the graphs of $g_u$ and $g_v$ respectively. 
By the choice of $U_4$ and  Corollary \ref{fonctionbornee2}, for all $v\in U_4$, $U_4=D(\zeta(v),v)$ for some $\zeta(v) \in br_e$. This means that for all $v\in U_4$, any element of the family $\{g_u|\Lambda_{\zeta(v)}; u\in U_4\setminus \{v\}\}$  agrees with $g_v$ only on $e$.  It follows by Lemma \ref{fonctionbornee} that there is some $\zeta\in br_e$ such that $\forall u, v\in U_4$, $u\neq v$, the functions $g_u|\Lambda_{\zeta}(e)$ and $g_v|\Lambda_{\zeta}(e)$ agree only on $e$. Replace $V_1$ by $V:=\Lambda_{\zeta}(e)$. \par\ 
  For every $x\in V^*$, the (definable) function from $U_4$ to $V$ which to $u$ associates $g_u(x)$ is now injective, thus a $C$-isomorphism on some neighborhood of $b$ (since all points in $V$ are independent of $b$) which we can suppose uniformly definable in $x$: 
  by $C$-minimality the union of all cones for which this is true, is a cone, 
  call it $J(x)$. 
For $\nu\in br_b$ define \[X(\nu):=\{x\in V:  \Lambda_{\nu}\subset J(x)\}.\]
By Corollary 2.4 there is a node $\nu_0$ on the branch of $b$ such that $X(\nu_0)=V$.
We replace $U_4$ by the cone of $b$ at $\nu_0$, which we call $U_5$. The family of functions $g_v:V\longrightarrow M$, $v\in U_5$ has the property $(*)$.\par\
Fix a generic element $u_0\in U_5$, and let $U_6$ be a subcone of $U_5$ containing $u_0$ such that $\forall u\in U_6$, $g_{u_0}(V)=g_u(V)$. For $u\in U_6$, let $f_u$ be the $C$-automorphism of $V$ defined by $f_u:=g_{u_0}^{-1}\circ g_u$. The family $\F:=\{f_u:u\in U_6\}$ is a nice family of functions.
  
 By $C$-minimality, there is a neighborhood $U$ of $u_0$ such that, either for all $u\in U\setminus\{u_0\}$, the function  $f_u$ is non-dilating, or for all $u\in U\setminus\{u_0\}$, the function  $f_u$ is dilating. 
Suppose the $f_u$ are dilating.  This means that 
$C_e((g_{u_0}^{-1}\circ g_u)^2,g_{u_0}^{-1}\circ g_u,id_V)$
 holds for any $u \in U$, $u \not= u_0$. 
 Define $D_u := \{ v \in U : C_e((g_{u}^{-1}\circ g_v)^2,g_{u}^{-1}\circ g_v,id_V) \}$ for $u \in U$. By genericity of $u_0$, for any $u$ close enough to $u_0$, $u \not= u_0$, $D_u$ contains $V_u \setminus \{u\}$ for some cone $V_u \ni u$.
Thus the function $U \rightarrow T(M) \setminus M$, $u \mapsto \inf \{ \mbox{basis  of }\Gamma ; \Gamma \mbox{ a cone},
\Gamma \setminus \{ u \}
\subset D_u, u \in \Gamma \}$, is well defined in the neighborhood of $u_0$. It must be locally constant.
So we can find two elements $u,v$ such that $u\in D_v$ and $v\in D_u$. It follows that \[C_e(g_{u}^{-1}\circ g_v,id_V, (g_{u}^{-1}\circ g_v)^{-1}) ,\] and \[C_e((g_{u}^{-1}\circ g_v)^{-1},id_V, g_{u}^{-1}\circ g_v) .\]  Contradiction.
 \hfill$\square$\par\ \par\noindent

\section{Tangents: existence and uniqueness}\label{stangents}
\begin{definition}Let $\F$ be a nice family of functions and $g:V\longrightarrow V$ a function such that $g(e)=e$. \begin{enumerate} \item Let $u$ be an element of $U$. The function $f_u$ is said to be \emph{tangent to $g$ relatively to the family $\F$} if for any $u'\in U$, we have $(\neg C)(f_u,f_{u'},g)$, in which case we write $f_u\sim_{\F} g$, or just $f_u\sim g$ if there is no confusion on $\F$.  
 \item The function $g$ is said to be \emph{derivable relatively to the family} $\F$ if there is a unique function $f_u\in \F$ such that $f_u\sim g$. In this case, $f_u$ is called \emph{the tangent to $g$ in $\F$}.
 \end{enumerate}
\end{definition}
{\bf\emph{Notation:}} We fix a nice family $\F$ of functions, and a definable function  $g:V \rightarrow V$  such that $g(e)=e$. We define \[T_g := \{ u \in U ; f_u\sim g  \},\] and for $u \in U$, \[\Gamma_{g,u} := \{ y \in U: C(f_u,f_y,g) \}.\]
\begin{lemma}\label{gammavide}
$\Gamma_{g,u} = \emptyset$ if and only if $u  \in T_g$. 
\end{lemma}
\emph{Proof:}
Fix elements $u$ and $v$ of $U$. So either $(\neg C)(f_u,f_v,g)$, or $C(f_u,f_v,g)$. By the definition of tangent, 
$u  \in T_g$ if and only if we are in the first case for every $v \in U$, if and only if $\Gamma_{g,u} = \emptyset$.\hfill$\square$\par\ \par\noindent
\begin{lemma}\label{Tforme}
\begin{enumerate}
\item Let $u$ and $v$ be elements of $T_g$. Then $\Lambda_{u \wedge v} \subset T_g$. Furthermore, there is a cone $W \subset V$ containing $e$ such that 
\[\forall z\in\Lambda_{u\wedge v}\ \forall x\in W: \Delta (f_u(x),f_z(x),g(x)).\] 
 \item Let $u \in T_g$ and $v \in U\setminus T_g$. Then the cone $\Gamma$ of $v$ at $u \wedge v$ is contained in $U\setminus T_g$. Furthermore, there is a cone $W \subset V$ containing $e$ such that \[\forall z\in\Gamma, \forall x\in \stern{W}:   C(f_z(x),f_u(x),g(x)).\]
\item If $T_g$ is not empty, then it is a ball. 
\end{enumerate}
\end{lemma}
\emph{Proof:}
\emph{1.}  Let $u,v$ be two distinct elements  of $T_g$, and let $W$ be a neighborhood of $e$ such that, \[\forall x\in W: \Delta(f_u(x), f_v(x), g(x)).\] Let $w\in \Lambda_{u \wedge v}$, so we have $\neg C(w,u,v)$. By $(*)$, we have that \[\forall x\in W:\neg C(f_w(x),f_u(x),f_v(x)).\]  So $\neg C(f_w(x), f_u(x), g(x))$ holds for all $x\in W$, and it follows that $w \in T_g$. For the second assertion, let $z$ be an element of $\Lambda_{u \wedge v}\subset T_g$. Let 
\[W_{z} := \bigcup \{  W: W \,\mbox {is a cone, }  
        e \in W\subset V, 
        \forall x \in W: 
        \Delta (f_u(x),f_z(x),g(x)) 
\}.\] $W_z$ is a non empty union of nested cones, so by $C$-minimality it is a cone. Let $\nu_z$ be its basis.
If $z' \in U$ is such that $C(u,z',z)$, then ($z' \in T_g$ and) $\nu_{z} = \nu_{z'}$. Thus the application  $z \mapsto \nu_{z}$ induces an application from the set of cones at $u\wedge v$ to the branch of $e$. As this set of cones equipped with the structure induced by $\M$ is strongly minimal and $br_e$ linearly ordered, the image of this application is finite. Let $\nu$ be the maximal element of the image. So we have 
\[\forall z\in \Lambda_{u\wedge v}, \forall x \in\Lambda_\nu: \Delta (f_u(x),f_{z}(x),g(x)).\] 
\emph{2.} Similar proof.\par\ \par\noindent
\emph{3.} By \emph{1}, either $T_g$ is empty, or it is a union of nested closed balls. It is then a cone or a closed ball by $C$-minimality. 
 \hfill$\square$\par\ \par\noindent
\begin{lemma}\label{chaine}
 Let $u,u',v,v'$ be elements of $U$ such that $v\in\Gamma_{g,u}$ and $v'\in\Gamma_{g,u'}$.  If $v$ is not an element of $\Gamma_{g,u'}$, then $v'$ is an element of $\Gamma_{g,u}$.
\end{lemma}
\emph{Proof:} Let $u,v,u',v'$ be elements of $U$ such that $v\in\Gamma_{g,u}$,  $v'\in\Gamma_{g,u'}$ and $v\notin\Gamma_{g,u'}$. Let $W\subset V$ be a neighborhood of $e$ such that, for every $x\in\stern{W}$, we have $C(f_u(x), f_v(x),g(x))$,  $C(f_{u'}(x), f_{v'}(x),g(x))$ and $\neg C(f_{u'}(x), f_{v}(x),g(x))$. By the first and third relations, we have \[\forall x\in\stern{W}: C(f_u(x),f_{u'}(x),g(x)).\] This together with the second relation yields \[\forall x\in\stern{W}: C(f_u(x),f_{v'}(x),g(x)).\]
 So  $v'$ is an element of $\Gamma_{g,u}$.
\hfill$\square$\par\ \par\noindent
\begin{lemma}\label{gammaforme}
If  $\Gamma_{g,u}$ is not empty, then it is a cone at a node on the branch of $u$. For a fixed $g$, the non empty $\Gamma_{g,u}$ form a chain of cones. 
\end{lemma}
\emph{Proof:} If $x$ and $y$ are elements of $\Gamma_{g,u}$, then $\Lambda_{x \wedge y}\subset \Gamma_{g,u}$.
So $\Gamma_{g,u}$ is a union of nested closed balls and, by $C$-minimality, it is a ball.
Fix an element $v\in \Gamma_{g,u}$. It is easy to see that $\Gamma_{g,u}$ contains the cone of $v$ at $u\wedge v$, so $\Gamma_{g,u}$ is a ball at a node $\nu$ on the branch of $u$. But it is clear that $u\notin\Gamma_{g,u}$. Thus $\nu=u\wedge v$ and  $\Gamma_{g,u}$ is the cone of $v$ at the node $u\wedge v$.  \par\noindent  The claim that the non empty $\Gamma_{g,u}$ form a chain of cones follows directly from Lemma \ref{chaine} and the fact that in $C$-structures, two cones have a nonempty intersection if and only if one of these cones is contained in the other one.\hfill$\square$\par\ \par\noindent
\begin{lemma}
There is an element $u\in U$ such that $f_u\sim g$. 
\end{lemma}
\emph{Proof:}
If $T_g$ is empty, then by Lemma \ref{gammavide}, no cone $\Gamma_{g,u}$ is empty. By Lemma \ref{gammaforme}, the $\Gamma_{g,u}$ form a chain of cones of $U$. By definable maximality, their intersection is not empty. But this intersection is contained in $T_g$, which is a contradiction. 
\hfill$\square$\par\ \par\noindent
\begin{lemma}\label{Wuniforme}
\begin{enumerate}
\item There is a cone $W$ containing $e$ such that $\Delta (f_u(x),f_v(x),g(x))$ holds for all $u,v \in T_g$ and all $x \in W$. 
\item Suppose that $T_g$ contains more than one element, and let $u \in T_g$. Then there is a cone $W$ containing $e$ such that $C(f_z(x),f_u(x),g(x))$ holds for all $z \in U \setminus T_g$ and all $x \in \stern{W}$. 
\end{enumerate}
\end{lemma}
\emph{Proof:}
\emph{(i).} By $C$-minimality and the definition of $T_g$, for all $u,v\in T_g$  there is some cone $W_{u,v}$ containing $e$ such that $\Delta (f_u(x),f_v(x),g(x))$ holds for all $x \in W_{u,v}$. Fixing $u$ and applying Lemma \ref{fonctionbornee} to the formula
$\varphi (i,v) = ``\Delta (f_u(x),f_v(x),g(x))\mbox{\it{ holds for all }} x \in \Lambda_i(e)$''
gives for all $u$, a cone $W_u$ containing $e$ such that $\Delta (f_u(x),f_v(x),g(x))$ holds for all $v \in T_g$ and all $x \in W_u$. Applying Lemma \ref{fonctionbornee} again gives the wanted result.  \par\ \noindent
\emph{(ii).} Suppose that $T_g$ is a ball containing at least two elements, and let $c_0 := \inf T_g$. Fix an element $a\in T_g$ and a cone $W$ containing $e$ such that $\Delta (f_a(x),f_z(x),g(x))$ holds for all $z \in T_g$ and all $x \in W$. For any node $c\leq c_0$ and any element $z \in U \setminus T_g$ such that $z\wedge a=c$,  there is a maximal ball $W_{c,z}\subset W$ such that  \[\forall x\in \stern{W}_{c,z}: C(f_z(x),f_a(x),g(x)).\] It is also easy to check that \[\forall x\in \stern{W}_{c,z}: C(f_{z'}(x),f_a(x),g(x))\] for every $z' \in U$ such that $C(a,z,z')$. Since $W_{c,z}\subset W$, we have 
\[\forall v\in T_g, \forall x\in \stern{W}_{c,z}: \Delta (f_{v}(x),f_a(x),g(x)) \wedge C(f_{z'}(x),f_a(x),g(x))\] 
for every $z'$ such that $C(a,z,z')$. Let $\nu_{c,z}$ be the basis of $W_{c,z}$. So the application  $z \mapsto \nu_{c,z}$ induces an application from the set of cones at $c$ not containing $T_g$ to the branch of $e$. By strong minimality of this set of cones, the image of 
this application is finite. Let $\nu_c$ be the maximal element of the image.  Now the application which to $c$ associates $\nu_c$ is an application from a bounded interval of the tree (namely the interval delimited by the basis of $U$ and that of $T_g$) to the branch of $e$. Its image is then bounded from above by some node $d$. If $W_0\subset W$ is the cone of $e$ at $d$, then we have  \[\forall u,v\in T_g,\forall z\in U\setminus T_g, \forall x\in \stern{W}_0: \Delta (f_{v}(x),f_u(x),g(x)) \wedge C(f_z(x),f_u(x),g(x)).\]
\hfill$\square$\par\ \par\noindent
\begin{prop}\label{unicite}
 The following are equivalent: \begin{enumerate}
\item there are elements $u,v \in U$ such that $C(f_u,f_v,g)$;
\item $T_g \not= U$; 
\item $g$ is derivable relatively to the family $\F$. 
\end{enumerate}
\end{prop}
\emph{Proof:} 
We show that (iii) follows from (i), the rest (namely (iii) $\Rightarrow$ (ii) $\Rightarrow$ (i)) is trivial. Let $u,v$ be like in the statement of (i). We know that $T_g$ is not empty. Suppose that it contains two distinct elements. By Lemma \ref{Tforme}(iii), $T_g$ is a cone or a closed ball. By Lemma \ref{Wuniforme}, there is a cone $W$ containing $e$ such that \[\forall \alpha\in T_g,\forall z\in U, \forall x\in W: \neg C(f_{\alpha}(x),f_z(x),g(x)).\,\,\,\,\,\,\,\,\,\,\,\,(a)\] 
Restricting $W$ if necessary, we can suppose that \[\forall x\in\stern{W}: C(f_u(x),f_v(x),g(x)).\]
Fix an element $x_0 \in \stern{W} $. By $(*)$, there is an element $w \in U$ such that $f_w(x_0)=g(x_0)$. Since $T_g$ contains more than one element, choose an element $\alpha\neq w$,  $\alpha\in T_g$. So we have $C(f_{\alpha}(x_0),f_w(x_0),g(x_0))$. This contradicts (\emph{a}).\hfill$\square$\par\ \par\noindent

\section{Derivability relative to nice families of functions}\label{sderivatives}
We fix a nice family $\F$ of functions.
\begin{lemma}\label{deriv2} Let $g:V\longrightarrow V$ be a function such that $g(e)=e$. Suppose that $g$ is derivable relatively to the family $\F$, and let $f_{u_1}$ be its tangent. Then for every $u\in U, u\neq u_1$, we have $C(f_u,f_{u_1},g)$.
\end{lemma}
\emph{Proof:} Follows directly from Lemma \ref{CCC} and the definition of derivability.
 \hfill$\square$\par\ \par\noindent
For $x\in \stern{V}$, the function $\phi_{\F,x}: U\rightarrow V, u\mapsto f_u(x)$, is a $C$-isomorphism from $U$ onto some cone $\phi_{\F,x}(U)$.\par\ \par\noindent
Let $g:V\longrightarrow V$ be a function such that $g(e)=e$. We define the map $\displaystyle\psi_{\F,g}: \stern{V}\rightarrow U, x\mapsto \phi_{\F,x}^{-1}(g(x))$.
So $\psi_{\F,g}(x)$  is the unique element $u$ of $U$ such that $f_u(x)=g(x)$, when such an element exists.
\begin{lemma}\label{psidefined}
 If $g$ is derivable relatively to the family $\F$, then $\psi_{\F,g}$ is defined on $\stern{W}$ for some neighborhood $W \subset V$ of $e$.
\end{lemma}
\emph{Proof:}
 Let $f_{u_1}$ be the tangent to $g$ in $\F$, and let $u\neq u_1$ be an element of $U$. By Lemma \ref{deriv2}, we have  $C(f_u,f_{u_1},g)$. So for some neighborhood $W \subset V$ of $e$, we have \[\forall x\in \stern{W}: C(f_u(x),f_{u_1}(x),g(x)).\] For every $x\in W$, $f_u(x)$ and $f_{u_1}(x)$ are elements of the cone  $\phi_{\F,x}(U)$, so the same holds for $g(x)$. Hence  $\psi_{\F,g}$ is defined on $\stern{W}$.
\hfill$\square$\par\ \par\noindent

\begin{prop}\label{limit} If $g$ is derivable relatively to $\F$, then $f_{u_1}\sim g$ if and only if \[\displaystyle \lim_{x\to e}\psi_{\F,g}(x)=u_1. \]
 \end{prop}
\emph{Proof:} By Lemma \ref{psidefined} $\psi_{\F,g}$ is defined on $\stern{W}_0$ for some neighborhood $W_0$ of $e$. For every $x\in W_0$, we have that $f_{\psi_{\F,g}(x)}(x)=g(x)$. Fix an element $u\in U\setminus\{u_1\}$ and a neighborhood $W_1 \subset V$  of $e$ such that \[\forall x\in \stern{W}_1: C(f_u(x), f_{u_1}(x),g(x)).\,\,\,\,\,\,\,\,\,\,\,\,\,\,\,\,(1)\]
From $(1)$ and the property $(*)$ of $\F$, it follows that \[C(u,u_1,\psi_{\F,g}(x)).\]
 This shows that the unique possible limit of $\psi_{\F,g}$ at $e$ is $u_1$. On the other hand, we know by definable maximality and Proposition $4.4$ of \cite{annalestoulouse}  that $\psi_{\F,g}$ has a limit at $e$. So we have that \[\displaystyle \lim_{x\to e}\psi_{\F,g}(x)=u_1. \]
\hfill$\square$\par\ \par\noindent

 \begin{definition}
  Let $\F$ and $\G$ be two nice families of functions having same domain and absorbing element.  
\begin{enumerate}\item The family $\G$ is said to be \emph{derivable relatively to} $\F$ if every element of $\G$  is derivable relatively to $\F$.
\item The families $\F$ and $\G$ are said to be \emph{comparable} if both are derivable relatively to each other.
\end{enumerate}
 \end{definition}
\noindent{\bf\emph{Notation:}} Let  $\F=\{f_u:u\in U\}$ and $\G\:=\{g_u:u\in U'\}$ be two nice families of functions. Suppose that the family $\G$ is derivable relatively to the family $\F$. We define the \emph{derivative} and denote it by $\partial_{\F,\G}$ as the function 

\begin{center}\begin{tabular}{ccll} $\partial_{\F,\G}:$&$U'$ &$\longrightarrow$ &$U$ \\ &$u'$  &$\longmapsto$  &$u\mbox{ such that } f_u\sim g_{u'}$.\end{tabular}
\end{center}
\begin{lemma}\label{deriveeiso1}
 Let  $\F=\{f_u:u\in U\}$ and $\G\:=\{g_u:u\in U'\}$ be two nice families of functions. Suppose that the family $\G$ is derivable relatively to the family $\F$. Let $a'_1, a'_2$ and $a'_3$ be elements of $U'$ such that $C(a'_1,a'_2,a'_3)$, and let $a_i:=\partial_{\F,\G}(a'_i)$.  Then 
 either $C(a_1,a_2,a_3)$ or $a_1=a_3$.
\end{lemma}
\emph{Proof:}
Suppose for a contradiction that $\neg C(a_1,a_2,a_3)$ and $a_3\neq a_1$ hold. Then $a_3\neq a_2$. Let $W \subset V$ be a neighborhood of $e$ such that, for every $x\in \stern{W}$ we have \[C(f_{a_3}(x), g_{a'_1}(x), f_{a_1}(x))\,\,\wedge\,\, C(f_{a_3}(x), g_{a'_2}(x), f_{a_2}(x))\] and \[ C(f_{a_1}(x), g_{a'_3}(x), f_{a_3}(x))\,\,\wedge\,\, C(f_{a_2}(x), g_{a'_3}(x), f_{a_3}(x)).\] Therefore, for every $x\in W$ holds: \[\neg C (g_{a'_1}(x), g_{a'_2}(x), g_{a'_3}(x)).\] By $(*)$ for $\G$ we have: \[\neg C(a'_1, a'_2, a'_3),\] contradiction.
\hfill$\square$\par\ \par\noindent
\begin{corollary}\label{deriveecontinue}
Let  $\F=\{f_u:u\in U\}$ and $\G\:=\{g_u:u\in U'\}$ be two nice families of functions. Suppose that the family $\G$ is derivable relatively to the family $\F$. Then $\partial_{\F,\G}$ is continuous on $U'$.
\end{corollary}
\emph{Proof:} By definable maximality, Proposition $4.4$ of \cite{annalestoulouse} and Lemma \ref{deriveeiso1}, the function $\partial_{\F,\G}$ admits a limit in each point of $U'$, and the limit is an element of $U$. Let $a'\in U'$ and  $a:=\partial_{\F,\G}(a')$. If the limit of $\partial_{\F,\G}$ in $a'$ is an element $b\neq a$, then we can find elements $b'_1,b'_2\in U'$ such that $C(b'_1,b'_2,a')$, and $C(a,b_i,b)$ for $i=1,2$, where $b_i=\partial_{\F,\G}(b'_i)$. Then $C(a,b_1,b_2)$, which contradicts Lemma \ref{deriveeiso1}.\hfill$\square$\par\ \par\noindent
\begin{prop}\label{symetrie1}
 Let $\F=\{f_u:u\in U\}$ and $\G\:=\{g_u:u\in U'\}$ be two nice families of functions. Suppose moreover that $\G$ is derivable relatively to $\F$. Let $a, a'$ be elements of $U$ and $U'$ respectively. Then $f_a\sim g_{a'}$ if and only if for any neighborhood $A\subset U$ of $a$, any neighborhood $A'\subset U'$ of $a'$, and any neighborhood $E\subset V$ of $e$, there are elements $\alpha\in A, \alpha'\in A'$ and $y\in \stern{E}$ such that $f_\alpha(y)=g_{\alpha'}(y)$. 
\end{prop}
\emph{Proof:}
For the implication from left to right  let $a\in U$, $a'\in U'$ and suppose that $f_a\sim g_{a'}$.  By Lemma \ref{limit}, we have that \[\displaystyle \lim_{x\to e}\psi_{\F,g_{a'}}(x)=a.\] Let $A\subset U,\,A'\subset U'$ and $E\subset V$ be neighborhoods of $a, a'$ and $e$ respectively. We can find an element $y\in \stern{E}$ such that $\psi_{\F,g_{a'}}(y)\in A$. Set $\alpha:=\psi_{\F,g_{a'}}(y)$ and $\alpha':=a'$. We have then that $\alpha\in A$, $\alpha'\in A'$, $y\in \stern{E}$, and by the definition of $\psi$ we have 
 \[f_\alpha(y)=g_{\alpha'}(y).\] 
 For the converse, let $b\in U$ be such that $f_b\not\sim g_{a'}$, and let $a\neq b$ be an element of $U$  such that $f_a\sim g_{a'}$. Let $E_0$ be a neighborhood of $e$ such that 
 \[\forall x\in \stern{E}_0: C(f_b(x),  g_{a'}(x), f_a(x)).\,\,\,\,\,\,\,\,\,\,\,\,(1)\] 
 Let $\Gamma$ be the cone of $a$ at the node $a\wedge b$. The element $\partial_{\F,\G}(a')=a$ is an element of $\Gamma$, so by continuity of $\partial_{\F,\G}$, there is an element $a'_1\in U'\setminus \{a'\}$ such that $a_1:=\partial_{\F,\G}(a'_1)\in \Gamma.$ 
 Let $E_1\subset V$ be a neighborhood of $e$  such that 
\[\forall x\in \stern{E}_1: C(f_b(x), g_{a'_1}(x), f_{a_1}(x)).\,\,\,\,\,\,\,\,\,\,\,\,(2)\] Let $A'$ be the cone of $a'$ at the node $a'\wedge a'_1$, $B$ be the cone of $b$ at the node $a\wedge b$, and $E:=E_0\cap E_1$.  By the properties of $\F$, we have \[\forall \beta\in B, \forall x\in \stern{E}: C(f_a(x),f_{\beta}(x),f_b(x) ).\,\,\,\,\,\,\,\,\,\,\,\,(3)\] By $(1)$,  $(2)$ and the fact that $C(f_b(x),f_{a_1}(x),f_a(x))$ for any $x \not= e$, we have \[\forall x\in \stern{E}: C(f_b(x),g_{a'}(x),f_a(x))\wedge C(f_b(x),g_{a'_1}(x)
,f_a(x) ).\,\,\,\,\,\,\,\,\,\,\,\,(4)\] Furthermore, by the properties of $\G$ we have: \[\forall \alpha\in A', \forall x\in \stern{E}: C(g_{a'_1}(x), g_{\alpha}(x), g_{a'}(x)).\,\,\,\,\,\,\,\,\,\,\,\,(5)\] From $(4)$ and $(5)$ we have \[\forall \alpha\in A', \forall x\in \stern{E}: C(f_b(x),g_{\alpha}(x),f_a(x) ).\,\,\,\,\,\,\,\,\,\,\,\,(6)\]
 By $(3)$ and $(6)$, we have \[\forall\beta\in B, \forall\alpha\in A', \forall x\in \stern{E}: f_{\beta}(x)\neq g_{\alpha}(x),\] and $B, A'$ and $E$ are neighborhoods of $b,a'$ and $e$ respectively.  This completes the proof. 
\hfill$\square$\par\ \par\noindent
An immediate consequence of Proposition \ref{symetrie1} is the
\begin{corollary}\label{symetrie2}
 Let $\F=\{f_u:u\in U\}$ and $\G=\{g_u:u\in U'\}$ be two comparable families of functions, and let $a, a'$ be elements of $U$ and $U'$ respectively. Then $f_a\sim_{\F} g_{a'}$ if and only if $g_{a'}\sim_{\G} f_a$.
\end{corollary}
\begin{corollary}\label{deriveeiso2} Let $\F=\{f_u:u\in U\}$ and $\G=\{g_u:u\in U'\}$ be two comparable families of functions. Then  $\partial_{\F,\G}$ and $\partial_{\G,\F}$ are inverse of each other, and they define (continuous) $C$-isomorphisms between $U$ and $U'$. \end{corollary}
\emph{Proof:} By Corollary \ref{symetrie2}, the functions $\partial_{\F,\G}$ and $\partial_{\G,\F}$ are inverse of each other, thus define bijections between $U$ and $U'$. Continuity follows from Corollary \ref{deriveecontinue}, and the rest from Lemma \ref{deriveeiso1}. \hfill$\square$\par\ \par\noindent
\begin{example}
 Let $\mathcal{M}$ be an algebraically closed valued field with maximal ideal $\Upsilon$, $U=U':=1+\Upsilon$, $V:=\Upsilon$, $u_0=u'_0=1$, $e=0$,  $\F:=\{u.x;u\in U\}$ and $\G:=\{x+(u-1).x^2;u\in U'\}$. Then $\G$ is derivable relatively to $\F$, and $\partial_{\F,\G}$ is constant and sends every $u'$ to $1$. Now $\partial_{\G,\F}$ is defined only at the point $1$, and its image in this point is $1$.
\end{example}
\begin{prop}\label{transitivite} Let $\F=\{f_u:u\in U\}$ and $\G=\{g_u:u\in U'\}$ be nice comparable families of functions, and $h:V\longrightarrow V$ a function such that $h(e)=e$. Suppose  that $h$ is derivable relatively to $\G$. Then $h$ is derivable relatively to $\F$. More precisely, if $f_a\sim_{\F}g_b\sim_{\G}h$, then $f_a\sim_{\F}h$.
\end{prop}
\emph{Proof:}
 Let $a,b$ be elements of $U,U'$ respectively such that $f_{a}\sim_{\F} g_{b}$ and $g_{b}\sim_{\G} h$. Suppose first that $\neg C(h,f_a,g_b)$ holds. 
For any $a_1\neq a$, we have $C(f_{a_1},f_a, h)$, which implies $f_a\sim_{\F} h$. \par\noindent
Now suppose $C(h,f_a,g_b)$, and $\neg C(f_{a_1},f_a,h)$ for some $a_1\neq a$. For every $y\neq b$, we have $C(g_y,g_b,h)$, thus $C(g_y,f_{a_1},g_b)$, and $g_b\sim_{\G} f_{a_1}$. From Corollary \ref{symetrie2} follows $f_{a_1}\sim_{\F} g_b$, so $a_1=a$. Contradiction.
\hfill$\square$\par\ \par\noindent

\section{The group law}\label{sgrouplaw}
Given a nice family of functions $\F=\{f_u:u\in U\}$, and  a $C$-automorphism $g$ of $V$ with $g(e)=e$, then $\{g\circ f_u:u\in U\}$ and $\{f_u\circ g:u\in U\}$ are nice families of functions.
\begin{lemma}\label{law1}Let $\F=\{f_u:u\in U\}$ be a nice family of non-dilating functions with identity $f_{u_0}$. Let $g,h$ be two $C$-automorphisms of $V$ with $g(e)=h(e)=e$. Assume that there is an element $c\in U$ such that $C(f_c,g,f_{u_0})$ and $C(f_c,h,f_{u_0})$. Then $C(f_c,g\circ h, f_{u_0})$ holds, and $g\circ h$ is derivable relatively to the family $\F$. If $f_x$ is the derivative of $g\circ h$, then $C(c,x,u_0)$.
 \end{lemma}

 \emph{Proof:}
We have  $C(g\circ f_c,g\circ h,g)$ (a), $C(f_c\circ f_c,g\circ f_c, f_c)$ (b) and $\neg C(f_c\circ f_c,f_c,id_V)$ (c). From (b) and (c) follows $C(id_V,g\circ f_c,f_c)$ (d), which together with $C(f_c,g,id_V)$ yields $C(g\circ f_c,g,id_V)$. By (a) we have $C(g\circ f_c,g\circ h,id_V)$. This relation together with (d) yields $C(f_c, g\circ h,id_V)$. The derivability of $g\circ h$ relatively to the family $\F$ follows from Proposition \ref{unicite}, and it is clear that if $f_x$ is the derivative of $g\circ h$, then $C(c,x,u_0)$.
\hfill$\square$\par\ \par\noindent

We fix a nice family of non-dilating functions $\F':=\{f_u:u\in U'\}$ with identity. Let $c\neq u_0$ be an element of $U'$, and $U:=\{x\in U':C(c,x,u_0)\}$. Then the family $\F:=\{f_u:u\in U\}$ is a nice family of non-dilating functions with identity.

\begin{lemma}\label{famoinsun}
Let $a\in U$. Then $C(f_c,f^{-1}_a,id_V)$.
\end{lemma}
\emph{Proof:} The relation $C(f_c,f_a,id_V)$ holds, so by composing with $f_c$ we have $C(f_c\circ f_c, f_a\circ f_c, f_c)$. This together with the fact that $f_c$ is not dilating yields
\[C(id_V, f_a\circ f_c, f_c). \,\,\,\,\,\,\,\,\,\,\,  (*)\]
Now we have $C(f_c,f_a,id_V)$ and that $f_a$ is not dilating.
From this follows that \[C(f_c,f_a\circ f_a,id_V). \,\,\,\,\,\,\,\,\,\,\,  (**)\]
The relations  ($*$) and ($**$) give $C(f_a\circ f_c,f_a\circ f_a,id_V)$.
Composing with $f^{-1}_a$ and using $C(f_c,f_a,id_V)$ gives the wanted result. 
\hfill$\square$\par\ \par\noindent 

\begin{lemma}\label{comparable}
Let $a$ be an element of $U$. Then the nice families $\{f_a\circ f_u:u\in U\}$, $\{f^{-1}_a\circ f_u:u\in U\}$, $\{f_u\circ f_a:u\in U\}$ and $\{f_u\circ f^{-1}_a:u\in U\}$ are comparable with $\F$.
\end{lemma}
\emph{Proof:}
Let $u\in U$. By Lemma \ref{law1}, $f_a\circ f_u$ is derivable relatively to $\F'$, with derivative $f_x$ for some $x\in U'$. By the second part of  Lemma \ref{law1}, $x\in U$, thus  $f_a\circ f_u$ is derivable relatively to $\F$. This shows that the family $\{f_a\circ f_u:u\in U\}$ is derivable relatively to $\F$.\par\ 
The derivability of $f_a \circ f_u$ relatively to $\cal F$ implies also that $f_u$ is derivable relatively to $\{ f_a^{-1} \circ f_u ; u \in U \}$. Lemma \ref{famoinsun} yields $C(f_c,f_a^{-1},id)$, so the previous argument applies to $f_a^{-1}$ instead of $f_a$. This shows that $\cal F$, $\{ f_a \circ f_u ; u \in U \}$ and $\{ f_a^{-1} \circ f_u ; u \in U \}$ are comparable.
That $\cal F$, $\{ f_u \circ f_a ; u \in U \}$ and $\{ f_u \circ f_a^{-1} ; u \in U \}$
are comparable can be proved in a similar way.
\hfill$\square$\par\ \par\noindent

{\bf\emph{Notations:}} If $g$ is derivable relatively to $\F$, we denote by $\widehat{g}$ its tangent. 

\begin{definition} Let $\comp:\F\times \F\rightarrow \F$ be the operation defined by $f_u\comp f_v:=\widehat{f_u\circ f_v}$. 
 \end{definition}

By Lemma \ref{comparable}, $\comp$ is well defined. We will show now that $(\F,\comp,f_{u_0})$ is a group (which is clearly infinite). This group structure can be obviously definably transferred on $U$.

\begin{lemma} The operation \emph{$\comp$} is regular.
 \end{lemma}

\emph{Proof:}
 We show right regularity, left regularity can be done in a similar way. Suppose that $f_w=f_u\comp f_{a}=f_v\comp f_{a}$. So $f_w\sim f_u\circ f_{a}$, and $f_w\sim f_v\circ f_{a}$. Let $\G:=\{f_x\circ f^{-1}_a: x\in U\}$. So $f_w\circ f^ {-1}_a\sim_{\G} f_u$. By Corollary \ref{symetrie2} we have  $f_u\sim f_w\circ f^ {-1}_a$. The same argument yields $f_v\sim f_w\circ f^ {-1}_a$, and we have $f_u=f_v$.
\hfill$\square$\par\ \par\noindent

\begin{lemma} The operation \emph{$\comp$} is associative.
 \end{lemma}

\emph{Proof:} 
If $g$ is such that  $f_a\sim g$, then $f_a\circ f_v\sim_{\G}  g\circ f_v$, where $\G:=\{f_x\circ f_v: x\in U\}$. Furthermore, it follows from the definition of $\comp$ that $f_a\comp f_v\sim f_a\circ f_v$. So by Proposition \ref{transitivite}, $f_a\comp f_v\sim g\circ f_v$. The same is true when we compose with $g$ from the right: $f_v\comp f_a\sim f_v\circ g$. This means with our notations: \[\widehat{g}\,\,\comp f_v=\widehat{g\circ f_v}\] and \[f_v\comp \,\,\widehat{g}=\widehat{f_v\circ g}.\]
This implies the following:
 \[(f_a\comp f_b)\comp f_c=\widehat{(f_a\circ f_b)}\comp f_c=\reallywidehat{\mbox{$(f_a\circ f_b)\circ f_c$}}=\reallywidehat{\mbox{$f_a\circ f_b\circ f_c$}}\] and 
\[f_a\comp (f_b\comp f_c)=f_a\comp (\widehat{f_b\circ f_c})=\reallywidehat{\mbox{$f_a\circ( f_b\circ f_c)$}}=\reallywidehat{\mbox{$f_a\circ f_b\circ f_c$}}.\]
So $\comp$ is associative.
\hfill$\square$\par\ \par\noindent

\begin{lemma}
 Every $f_a\in\F$ admits an inverse.
\end{lemma}

\emph{Proof:}
 Let $\G$ be the family $\{f_a\circ f_u: u\in U\}$. The families $\F$ and $\G$ are comparable, so let $b\in U$ be such that $f_a\circ f_b\sim_{\G} f_{u_0}$. By Corollary \ref{symetrie2} we have $f_{u_0}\sim f_a\circ f_b$, so $f_a\comp f_b=f_{u_0}$, which is the identity element of $(\F,\comp)$. So $f_a$ has a right inverse, and the same argument shows that $f_a$ has a left inverse.\hfill$\square$\par\ \par\noindent

We have proved the following

\begin{prop}\label{famillegroupe}Let $\M$ be a definably maximal and non-trivial geometric $C$-minimal structure. Suppose moreover that in the underlying tree of $\M$ there is no definable bijection between a bounded interval and an unbounded one, and that  a nice family of non-dilating functions with identity is definable in $\M$. Then there is an infinite group definable in $\M$.
\end{prop}

 Propositions \ref{existencefamille} and \ref{famillegroupe} yield directly the following Theorem.

\begin{theorem}Let $\M$ be a definably maximal and non-trivial geometric $C$-minimal structure. Suppose moreover that in the underlying tree of $\M$ there is no definable bijection between a bounded interval and an unbounded one. Then there is an infinite group definable in $\M$.
 \end{theorem}

{\bf\emph{acknowledgements:}}
  We are indebted to Bernhard Elsner for pointing out an error in the original proof of Proposition \ref{existencefamille}.

 \nocite{adeleke}
\nocite{macphersonhaskell}
\nocite{macphersonsteinhorn}
\nocite{maalouf1}

\noindent
   Fran\c coise Delon, CNRS, \\
   IMJ-PRG, UMR 7586, \\Univ Paris Diderot,  F-75013, Paris, \\
   UFR de math\'ematiques, case 7012\\
75205 Paris Cedex 13 \\
France \\
   email: delon@math.univ-paris-diderot.fr
\par\ \par\noindent
   Fares Maalouf\\
  Universit\'e St Joseph\\
ESIB\\
B.P. 11-514\\ 
Riad el Solh, Beirut 11072050\\
Lebanon     \\
email:fares.maalouf@usj.edu.lb

\end{document}